\documentclass[reqno,10pt,a4paper]{amsart}
\usepackage{srcltx}
\usepackage{graphicx,color,amssymb,latexsym,amsmath,amsfonts}
\usepackage{mathrsfs}

\newtheorem{thm}{Theorem}[section]
\newtheorem{proposition}{Proposition}[section]
\newtheorem{lemma}{Lemma}[section]

\newtheorem{ex}{Example}[section]
\newtheorem{remark}{Remark}[section]
\numberwithin{equation}{section}

\input cyracc.def

\begin{document}

\title{On minimax nonparametric estimation \\ of signal
in Gaussian noise}

\
\author{Mikhail Ermakov, IPME RAS, RUSSIA}

\begin{abstract} For the problem of nonparametric estimation of signal in Gaussian noise we point out the strong asymptotically minimax estimators on maxisets for linear estimators (see \cite{ker93,rio}). It turns out  that the order of rates of convergence of Pinsker estimator on this maxisets is worse than the order of rates of convergence for the class of linear estimators considered on this maxisets. We show that balls in  Sobolev spaces are maxisets for Pinsker estimators.
\end{abstract}


\maketitle



\section{\bf  Introduction \label{s1}}
For nonparametric estimation problem of  a signal in Gaussian white noise  optimal rates of convergence of estimators has been explored for a wide range of functional spaces and for a completely different setups (see \cite{don, nem,jo,ts} and references therein). The strong asymptotically minimax estimators are known only if a priori information is provided that a signal belongs to ellipsoid in $L_2$ \cite{kuk,pin,jo,ts,nus}, balls in $L_\infty$ \cite{ber,don94,kor,lep} or for bodies in Besov spaces \cite{jo}.  The paper goal is to pay attention that strong asymptotically minimax estimators can be obtained for other sets of functions. For trigonometric orthogonal system of functions the definition of these sets coincides with the definition of a ball in Besov space $B^\alpha_{2\infty}$ for some norm. We shall denote these sets as $B(\alpha,P_0)$ with $\alpha>0$ and $P_0 > 0$.

The balls $B_{2\infty}^\alpha(P_0)$ have remarkable properties in nonparametric estimation. \vskip 0.3cm
This sets carry a rather reasonable information on a signal smoothness.
\vskip 0.3cm
These sets are the  sets having a given rates of convergence for the most wellknown linear nonparametric estimators\cite{ker93, ker02}.
\vskip 0.3cm
For linear statistical estimators these sets are the largest sets  with a given rate of convergence \cite{rio}.
\vskip 0.3cm
Nonparametric estimation of solutions of linear ill-posed inverse problems in Gaussian noise for the sets $B(\alpha,P_0)$ has been also explored earlier in econometrics \cite{rio09}.

The arising strong asymptotically minimax estimators are penalized maximum likelihood estimators for some quadratic penalty function \cite{eg},\cite{wah}. Thus we obtain that  likelihood estimation with quadratic penalty function is optimal not only in Bayes sense but  in the minimax sense as well. These asymptotically minimax estimators are also trigonometric spline estimators \cite{jo, ts, wah}. The results can be also interpreted as a solution of inverse problem. For Bayes estimators and  maximum penalized estimators one needs to find the largest sets such that these estimators are asymptotically minimax on these sets.

The nonasymptotic setup is also explored. In this setup we show that our estimator is minimax for the class of all linear estimators.

The results can be easily modified on the case of minimax estimation of linear ill-posed  problem. For this setup minimax estimator can be treated as some version of Tikhonov regularization algorithm \cite{tik}.

We show that the order of rates of convergence of Pinsker estimator on $B(\alpha,P_0)$ is worse than the order of rates of convergence for the asymptotically minimax estimators  on this maxisets, and the balls in  Sobolev spaces are maxisets for Pinsker estimators.

The  results will be provided in terms of sequence model. Let we observe  a random sequence $y= \{y_j\}_{j=1}^\infty$,
$$
y_j = x_j + \epsilon\,\sigma_j\xi_j, \quad \epsilon > 0,\quad 1 \le j < \infty
$$
where $\sigma_j>0$ are known constants and $\xi_j, 1 \le j < \infty,$ are independent Gaussian random variables, $E\xi_j = 0$ and $E \xi^2_j =1$.

The problem is to estimate the parameter $x = \{x_j\}_{j=1}^\infty$.

Denote $\sigma = \{\sigma_j\}_{j=1}^\infty$ and $\xi = \{\xi_j\}_{j=1}^\infty$.

For the estimation with  fixed $\epsilon >0$ minimax estimators will be established if a priori information is provided in the following form
\begin{equation}\label{r2}
x \in B(a,P_0) = \left\{x= \{x_i\}_{i=1}^\infty : \sup_k a_k^{-1} \sum_{j=k}^\infty x_j^2 \le P_0\right\}
\end{equation}
where $a=\{a_k\}_{k=1}^\infty$ and $a_k> 0$ is decreasing sequence.

For asymptotically minimax estimation we shall consider the more narrow class of sets $B(\alpha,P_0) = B(\tilde a,P_0)$ with $\tilde a=\{k^{-2\alpha}\}, \alpha>0$.
The analysis of the proof shows that the results can be extended on another sequences $a_k$.
However this requires more accurate reasoning. For trigonometric orthogonal system of function  $\sup_k a_k^{-1} \sum_{j=k}^\infty x_j^2$ can be considered as some norm in Besov space $B^\alpha_{2\infty}$. For Besov bodies in $B^r_{2\infty}$ generated wavelets asymptotically minimax estimators one can find in Johnstone \cite{jo}. For this setup another extremal problem arises.

There are numerous research on strong adaptive asymptotically minimax estimation \cite{ jo,ts}.
The results on adaptive estimation in Pinsker model \cite{jo,ts} are easily carried over on  paper setup for asymptotically minimax estimation on the sets $B(\alpha,P_0)$.

Below we remind the definition of maxisets.

For estimator $\hat x_\epsilon$, for the loss function $||\hat x_\epsilon - x ||^2$, for rates of convergence $\epsilon^\gamma, \gamma > 0$, and  for the constant $C>0$, the maxiset is
$$
 MS(\hat x_\epsilon,\gamma)(C) = \{ x : \sup_\epsilon \epsilon^{-2\gamma}  E_x ||\hat x_\epsilon -x ||^2 < C \}.
 $$
 Here $\|x\|$ denotes norm of vector $x=\{x_j\}_{j=1}^\infty$  in Hilbert space,  $$\|x\|^2 = \sum_{j=1}^\infty x_j^2.$$

In what follows we shall denote letters $c, C$ positive constants and let $ a_\epsilon \asymp b_\epsilon$ imply $c < a_\epsilon/b_\epsilon < C$ for all $\epsilon > 0$.
\section{\bf Main Results}
We say that linear estimator $\hat x_\epsilon = \{\hat x_{\epsilon j}\}_{j=1}^\infty$ is minimax  in the class of linear estimators $\hat x_{\epsilon\lambda} = \{\hat x_{\epsilon j\lambda_j}\}_{j=1}^\infty, \hat x_{\epsilon j\lambda_j} = \lambda_j y_j, \lambda_j \in R^1, 1 \le j < \infty, \lambda= \{\lambda_j\}_{j=1}^\infty$,  if
\begin{equation}\label{r3}
R_{l\epsilon}\doteq\sup_{x \in B} E_x||\hat x_\epsilon - x ||^2 =
\inf_{\lambda}\sup_{x \in B} E_x||\hat x_{\epsilon\lambda} - x ||^2.
\end{equation}
We say that the estimator $\hat x_\epsilon$ is asymptotically minimax if
\begin{equation}\label{}
R_\epsilon\doteq \sup_{x \in B(\alpha,P_0)} E_x||\hat x_\epsilon - x ||^2 =
\inf_{\tilde x_\epsilon \in \Psi}\sup_{x \in B^r_{2\infty}(P_0)} E_x||\tilde x_\epsilon - x ||^2(1 + o(1))
\end{equation}
as $\epsilon \to 0$. Here $\Psi$ is the set of all estimators.

The minimax estimator in the class of linear estimators will be established if the following assumptions hold.

\noindent{\bf A1} There is $c>0$ such that $c< \sigma_j^2 < \infty$ for all $j$.

\noindent{\bf A2}. For all $j>1$
\begin{equation}\label{u1}
\frac{\sigma_j^2 (a_{j-1}- a_j)}{\sigma_{j-1}^2 (a_j - a_{j+1})} > 1.
\end{equation}
This implies that sequence $\sigma_j^2 (a_{j-1}- a_j)$ is strictly increasing.
\begin{thm}\label{t1} Assume A1,A2. Then the linear estimator $\hat x_\lambda$  with
\begin{equation}\label{u2}
\lambda_j = \frac{P_0(a_j- a_{j+1})}{P_0(a_j- a_{j+1}) + \epsilon^2\sigma_j^2}.
\end{equation}
is  minimax on the set of all linear estimators.

The  minimax risk equals
\begin{equation}\label{u3}
R_{l\epsilon} = \epsilon^2\sum_{j=1}^\infty \frac{P_0\sigma_j^2(a_j- a_{j+1})}{P_0(a_j- a_{j+1}) +\epsilon^2\sigma_j^2 }.
\end{equation}
\end{thm}
\begin{remark} The estimator $\hat x_\lambda$ is maximum penalized likelihood estimator \cite{eg,wah} with quadratic  penalty function
$$
P_0^{-1} \sum_{j=1}^\infty(a_j - a_{j+1})^{-1} \sigma_j^2 x_j^2
$$
and Bayes estimator  with a priori measure corresponding independent Gaussian coordinates $x_j$, $E x_j = 0$ and $E x_j^2 = P_0(a_j - a_{j+1})$, $1 \le j < \infty$.
\end{remark}
\begin{remark} Theorem \ref{t1} holds also for a finite number of  observations $y_1,\ldots,y_k$ with $k < \infty$.
\end{remark}
In Theorem \ref{t2} we replace A2 more simple assumption.

\noindent{\bf B1}. For all $j>j_0$
\begin{equation}\label{}
\frac{\sigma_j^2j^{2\alpha+1}}{\sigma_{j-1}^2 (j-1)^{2\alpha+1}} > 1.
\end{equation}
This implies that sequence $\sigma_j^2j^{2\alpha+1}$ is strictly increasing.
\begin{thm}\label{t2} Assume A1,B1. Then the linear estimator  $\hat x_\lambda$ with
\begin{equation}\label{}
\lambda_j = \frac{2\alpha P_0j^{-2\alpha-1}}{2\alpha P_0j^{-2\alpha-1} + \epsilon^2\sigma_j^2}.
\end{equation}
is asymptotically minimax on the set of all estimators.

The asymptotically minimax risk equals
\begin{equation}\label{}
R_\epsilon = \epsilon^2\sum_{j=1}^\infty \frac{2\alpha P_0j^{-2\alpha-1}\sigma_j^2}{2\alpha P_0j^{-2\alpha-1} + \epsilon^2\sigma_j^2 }(1+o(1)).
\end{equation}
\end{thm}

\begin{remark} The estimator $\hat x_\lambda$ is maximum penalized likelihood estimator \cite{eg,wah} with quadratic  penalty function
$$
(2\alpha P_0)^{-1} \sum_{j=1}^\infty j^{1+2\alpha} \sigma_j^2 x_j^2
$$
and Bayes estimator  with a priori measure corresponding independent Gaussian coordinates $x_j$,  $E x_j = 0$ and $E x_j^2 = 2\alpha P_0 j^{-1-2\alpha}$, $1 \le j < \infty$.
\end{remark}
 Theorems \ref{t1} and \ref{t2} are extended easily on linear ill-posed inverse problem setup. The maxisets for for linear ill-posed inverse problems has been studied Loubes  and Rivoirard \cite{rio09}.

  Suppose we observe a random vector
$$ y = Rx +\epsilon\xi$$ with linear self-adjoint operator $R: H \to H$ in a separable Hilbert space $H$. Other notations are the same as in previous setup.

Suppose the linear operator $R$ admits singular value decomposition (see \cite{ts,jo, ing,rio09}) with eigenvalues $r_j, 1 \le j < \infty$.  Then we can consider this setup in the following form.

We observe random vector
$$ z_j = r_j x_j + \epsilon\sigma_j\xi_j, \quad 1 \le j <\infty. $$ Suppose $\xi_j, 1 \le j <\infty,$ are i.i.d. Gaussian r.v.'s, $E\xi_j =0, E\xi_j^2 =1$. The problems of estimation of $x = \{x_j\}_{j=1}^\infty$ are the same. Dividing on $r_j$, we obtain the setup of signal estimation.
\vskip 0.3cm

Below two asymptotics of minimax risks for linear ill-posed inverse problems are provided.
\begin{ex}\label{e1} Let $\alpha >0, \gamma > 0$. Let $|r_j| = C j^{-\gamma}(1 + o(1))$ and $\sigma_j =1, 1 \le j<\infty$. Then
\begin{equation}\label{}
R_\epsilon= \epsilon^{\frac{4\alpha}{1+2\alpha+2\gamma}}\frac{\pi}{2\alpha
\sin\left(\frac{\pi(2\gamma+1)}{2\alpha}\right)}(2\alpha P_0)^{\frac{2\gamma+1}
{2\gamma+2\alpha+1}}C^{-\frac{2\alpha}{2\gamma+2\alpha+1}}(1+o(1)).
\end{equation}
\end{ex}
\begin{ex}\label{e1} Let  $\alpha >0, \gamma >0, B>0, \kappa \in R^1$. Let $|r_j| = Cj^{-\kappa}\exp\{-Bj^\gamma\}$ and $\sigma_j =1, 1 \le j<\infty$. Then
\begin{equation}\label{}
R_\epsilon= P_0 B^{2\alpha/\gamma}|\log \epsilon|^{-2\alpha/\gamma}(1 + o(1)).
\end{equation}
\end{ex}
Note that these asymptotics coincide with the asymptotics of risks of corresponding Bayes estimators.

Johnstone  (Th 3.10, Ch3, \cite{jo}) has provided the comparison of strong asymptotics of minimax risks for trigonometric spline estimators  and Pinsker estimators if unknown signal  belongs to a ball in Sobolev space. The trigonometric spline estimators are strong asymptotically minimax estimators on maxisets $B(\alpha,P_0)$. Thus we can consider this result as a comparison of risk asymptotics for strong asymptotically minimax estimators on maxisets $B(\alpha,P_0)$ and Pinsker estimators. Below we provide similar comparison, if a priori information is provided, that unknown signal belongs to maxiset $B(\alpha,P_0)$.

Pinsker estimator $\tilde\theta_{\epsilon\mu} = \{\tilde\theta_{\epsilon j}\}_{j=1}^\infty$ is linear estimator
$$
\tilde\theta_{\epsilon j} = \lambda_{\epsilon j} y_j
$$
with
$$
\lambda_{\epsilon j} = (1 -\mu b_j)_+
$$
where $b_j =j^\beta, \beta>0,$ and parameter $\mu$ is defined by equation
$$
\epsilon^2 \sum_{j=1}^\infty b_j^2((\mu b_j)^{-1} -1)_+ = P.
$$
Pinsker estimator is asymptotically minimax on ellipsoids
$$
S(\beta,P) =\left\{x : \sum_{j=1}^\infty b_j^2 x_j^2 \le P, x =\{x_j\}_1^\infty\right\},
$$
with $b =\{b_j\}_1^\infty$ and $P>0$.

Denote
$$
R_\epsilon(\alpha,\beta) = \inf_\mu\sup_{\theta\in B(\alpha,P_0)} E_\theta  ||\tilde\theta_{\epsilon\mu} -\theta||^2.
$$
Denote $$ C = \frac{2\alpha^2}{(1+\alpha)(1+2\alpha)}$$
\begin{thm}\label{t3} Let $0<\alpha < \beta$. Then
\begin{equation}\label{}
R_\epsilon(\alpha,\beta) = C^{\frac{2\alpha}{1+2\alpha}} C_1^{\frac{1}{1+2\alpha}}\left((2\alpha)^{-\frac{2\alpha}{1+2\alpha}} + (2\alpha)^{\frac{1}{1+2\alpha}}\right)\epsilon^{\frac{4\alpha}{1+2\alpha}}
\end{equation}
with $C_1 = \frac{\beta}{\beta-\alpha}P_0$.

Let $\alpha > \beta>0$. Then
\begin{equation}\label{}
R_\epsilon(\alpha,\beta) = C^{\frac{2\beta}{1+2\beta}} C_1^{\frac{1}{1+2\beta}}\left((2\beta)^{-\frac{2\beta}{1+2\beta}} + (2\beta)^{\frac{1}{1+2\beta}}\right)\epsilon^{\frac{4\beta}{1+2\beta}}
\end{equation}
with  $$C_1 = \sum_{j=1}^\infty j^{2\beta}\left(j^{-2\alpha}  - (j+1)^{-2\alpha}\right)$$

If $\alpha = \beta$, then
\begin{equation}\label{}
R_\epsilon(\alpha,\beta) = \left((2\alpha^2)^{\frac{1}{1+2\alpha}} + 2^{-\frac{2\alpha}{1 + 2\alpha}}\alpha^{\frac{1-2\alpha}{1+2\alpha}}\right) (1+2\alpha)^{-\frac{1}{1+2\alpha}} P_0^{\frac{1}{1+2\alpha}}C^{\frac{2\alpha}{1+2\alpha}}\epsilon^{\frac{4\alpha}{1+2\alpha}} |2\ln \epsilon|^{\frac{1}{1+2\alpha}}.
\end{equation}
\end{thm}
The most interest represents the comparison of risks of Pinsker estimator and asymptotically minimax estimators on maxisets if $\alpha= \beta$. For this setup we compare the risks of estimators on the sets having almost the same smoothness. We see that the risks of Pinsker estimators have additional logarithmic term in asymptotic. Pinsker estimators do not belong to the class of linear estimators having the maxisets $B(\alpha,P_0)$. It turns out that the balls in Sobolev space $S(\beta,P)$ are maxisets for Pinsker estimators.
\begin{thm}\label{t4} There exists $C>0$ such that, for all $\epsilon>0$,
\begin{equation}\label{}
R_\epsilon(\beta,x) = \epsilon^{-\frac{4\beta}{1 +2\beta}}\inf_\mu E_x ||\tilde x_{\epsilon\mu} - x ||^2 < C < \infty,
\end{equation}
if and only if $x$ belongs to Sobolev space
$$
S^\beta = \left\{x : \sum_{j=1}^\infty b_j^2x_j^2 < \infty, \quad x = \{x_j\}_{j=1}^\infty\right\}.
$$
\end{thm}
In the theory of linear ill-posed inverse problems one of the most wide spread assumption is that the solution  $x$ satisfies a source condition \cite{cav}
$$
x \in \{x:\, x = B u,\, \| u \| \le 1,\, u \in H \},
$$
where $B$ is linear self-conjugate compact operator. This implies that the solution   $x$ belongs to ellipsoid.  Theorems \ref{t3}  and \ref{t4} show that optimal linear solution on such sets can have worse rates of convergence on more wider sets then other linear estimators. 
\section{Proof of Theorems}
\subsection{\bf Proof of Theorem \ref{t1}}
We begin with the proof of lower bound. Denote $\theta_j^2 = P_0(a_j- a_{j+1}), \theta= \{\theta_j\}_{j=1}^\infty$.

We have
\begin{equation}\label{}
\inf_\lambda\sup_{x\in B} E_x || \hat x_\lambda- x||^2 \ge \inf_\lambda E_\theta ||\hat x_{\epsilon\lambda} - \theta||^2\\=
\epsilon^2\sum_{j=1}^\infty\frac{\theta_j^2\sigma_j^2}{\theta_j^2 + \epsilon^2\sigma_j^2}
\end{equation}
and infimum is attained for
$$
\lambda_j = \frac{\theta_j^2}{\theta_j^2 + \epsilon^2\sigma_j^2} = \frac{P_0\,(a_j- a_{j+1})}{P_0\,(a_j- a_{j+1}) + \epsilon^2\,\sigma_j^2}.
$$
Proof of upper bound is based on the following reasoning. Let $x=\{x_j\}_{j=1}^\infty \in B$. For all $k$ denote
$$ u_k = a_k^{-1}\sum_{j=k}^\infty x_j^2.$$
Then $x_k^2= a_k u_k - a_{k+1} u_{k+1}$.

For the sequence of $\lambda_j$ defined in Theorem \ref{t1}, we have
\begin{equation}\label{ux1}
\begin{split}&
E_x\sum_{j=1}^\infty\,(\lambda_jy_j - x_j)^2 =\epsilon^2\,\sum_{j=1}^\infty \lambda_j^2\sigma_j^2 + \sum_{j=1}^\infty(1 -\lambda_j)^2x_j^2\\&= \epsilon^2 \sum_{j=1}^\infty \lambda_j^2\sigma_j^2 +  \sum_{j=1}^\infty(\theta_j^2\sigma_j^{-2}\epsilon^{-2} +1)^{-2}(a_j u_j - a_{j+1} u_{j+1})\\&=\epsilon^2 \sum_{j=1}^\infty \lambda_j^2\sigma_j^2  + (\theta_1^2\sigma_1^{-2}\epsilon^{-2} +1)^{-2} u_1
\\&-
\sum_{j=2}^\infty u_j a_j\left((\theta_{j-1}^2\sigma_{j-1}^{-2}\epsilon^{-2} +1)^{-2} -
(\theta_j^2\sigma_j^{-2}\epsilon^{-2} +1)^{-2}\right).
\end{split}
\end{equation}
By A2, the last addendums in the right hand-side of (\ref{ux1}) are negative. Therefore the supremum of right hand-side of (\ref{ux1}) is attained for $u_j= P_0$, $1 \le j< \infty$. This completes the proof of Theorem  \ref{t1}.
\subsection{\bf Proof of Theorem \ref{t2}}
The upper bound follows from Theorem \ref{t1}. Below the proof of lower bound will be provided. This proof has a lot of common features with the proof of lower bound in Pinsker Theorem  \cite{jo, pin, ts}.

Fix values $\delta_1, 0<\delta_1<1,$ and $\delta, 0 < \delta< P_0$. Define a family of natural  numbers $k_\epsilon, \epsilon>0,$ such that $\epsilon^{-2} \sigma^2_{k_\epsilon} 2rP_0k_\epsilon^{-2r-1} = 1 + o(1)$ as $\epsilon \to 0$. Define sequence $\eta = \{\eta_j\}_{j=1}^\infty$ of Gaussian i.i.d.r.v.'s $\eta_j = \eta_{j\delta\delta_1}, E[\eta_{j}] = 0, \mbox{Var}[\eta_j] = (P_0 - \delta)(2r)^{-1} j^{-2r-1}$, if $\delta k_\epsilon \le j \le \delta^{-1} k_\epsilon$, and $\eta_j=0$ if $j < \delta_1 k_\epsilon$ or $j > \delta_1^{-1}k_\epsilon$.

Denote $\mu$ the probability measure of random vector $\eta$. Define $\tilde x$ Bayes estimator with a priory measure $\mu$.

 Define the conditional probability measure $\nu_\delta$ of random vector
  $\eta$ given $\eta \in  B(\alpha,P_0).$
  Define $\bar x $ Bayes estimator of $x$ with a priori measure $\nu_\delta$.
 Denote $\theta$ the random variable having probability measure $\nu_\delta$.

For any estimator $\hat x$ we have
\begin{equation}\label{e1}
\begin{split}&
\sup_{x \in B^s_{2\infty}} E_x||\hat x -x||^2 \ge E_{\nu_\delta}E_\theta ||\hat x - \theta||^2\\&\ge
E_\mu E_\eta ||\tilde x -\eta||^2 -E_\mu E_\eta(||\bar x - \eta||^2, \eta \notin B(\alpha,P_0)) P^{-1}_\mu(\eta \in B(\alpha,P_0)) .
\end{split}
\end{equation}
We have
\begin{equation}\label{e2}
 {\mathbf E}_\mu\,{\mathbf E}_\eta\,\|\tilde x -\eta\|^2 = I(P_0-\delta)(1+o(1)),
\end{equation}
where
$$
I(P_0-\delta) =
\epsilon^2\sum_{j=l_1}^{l_2}\,
\frac{\sigma_j^2}{1 + (2\alpha\,(P_0-\delta_1))^{-1}\epsilon^2\,\sigma_j^2\,j^{2\alpha+1}}
$$ 
with $l_1 = [\delta_1 k_\epsilon]$ and $l_2 =[\delta_1^{-1}k_\epsilon]$. Here $[a]$ denotes whole part of a number $a \in R^1$.

Since
 $$||\bar x ||^2 \le \sup_{x\in B^r_{2\infty}} ||x||^2 \le P_0,$$
 we have
\begin{equation}\label{e3}
\begin{split}&
E_\mu E_\eta(||\bar x - \eta||^2, \eta \notin B(\alpha,P_0)) \le  2 E_\mu E_\eta ( ||\bar x||^2 + ||\eta||^2, \eta \notin B(\alpha,P_0)) \\& \le 2 P_0 P_\mu(\eta \notin B(\alpha,P_0)) + \sum_{j=l_1}^{l_2} (E_\mu \eta_j^4)^{1/2} P_\mu^{1/2}(\eta \notin B(\alpha,P_0)).
\end{split}
\end{equation}
Since
$E_\mu [\eta_j^4] \le C j^{-2(\alpha-2},$ we have
\begin{equation}\label{e4}
\sum_{j=l_1}^{l_2} (E_\mu \eta_j^4)^{1/2} \le C\delta_1^{-\alpha} k_\epsilon^{-2\alpha}.
\end{equation}
It remains to estimate
\begin{equation}\label{e5}
P_\mu(\eta \notin B^r_{2\infty}) = P\left(\max_{l_1 \le i \le l_2} i^{2\alpha} \sum_{j = i}^{l_2} \eta_j^2-P_0(1-\delta_1/2) > P_0\delta_1/2\right) \le \sum_{i=l_1}^{l_2} J_i
\end{equation}
with
$$
J_i = P\left( i^{2\alpha} \sum_{j = i}^{l_2} \eta_j^2-P_0(1-\delta/2)> P_0\delta/2\right)
$$
To estimate $J_i$ we implement the following Proposition \cite{hs}.
\begin{proposition}\label{p1} Let $\xi = \{\xi_i\}_{i=1}^l$ be Gaussian random vector with i.i.d.r.v.'s $\xi_i$, $E \xi_i = 0, E \xi^2 =1$. Let $A$ be $l\times l$-- matrix and $\Sigma = A^T A$. Then
\begin{equation}\label{e6}
P(||A\xi||^2 > \mbox{tr}(\Sigma) + 2\sqrt{\mbox{tr}(\Sigma^2)t} + 2 ||\Sigma||t) \le \exp\{-t\}.
\end{equation}
Here $\mbox{\rm tr}\,(\Sigma)$ denote the trace of matrix $\Sigma$.
\end{proposition}
Define matrix $\Sigma= \{\sigma_{lj}\}_{l,j=i}^{l_2}$ with $\sigma_{jj} = j^{-2\alpha-1}i^{2\alpha}\frac{P_0-\delta}{2\alpha}$ and $\sigma_{lj} =0$ if $l\ne j$. Then
\begin{equation}\label{e7}
2\sqrt{\mbox{tr}(\Sigma^2)t} + 2 ||\Sigma||t =  \frac{P_0-\delta}{\alpha(4\alpha+1)} \sqrt{i^{-1}t}(1+o(1)) + i^{-1}t \doteq V_i(t)
\end{equation}
We put $t = k_\epsilon^{1/2}$. Then $V_i(t) < Ck_\epsilon^{-1/2}, l_1 \le i \le l_2$ and implementing (\ref{e6})  we have
\begin{equation}\label{e8}
J_i < \exp\{-k_\epsilon^{-1/2}\}
\end{equation}
and therefore
\begin{equation}\label{e8a}
\sum_{j=l_1}^{l_2} J_i \le \delta_1^{-1}k_\epsilon\exp\{- k_\epsilon^{1/2}\}
\end{equation}
To complete the proof it remains to estimate $R_\epsilon-I(P_0-\delta)$.

By straightforward estimation, it is easy to verify that
\begin{equation}\label{e9}
|I(P_0) - I(P_0 -\delta)| < C\delta I(P_0)
\end{equation}
We have
\begin{equation}\label{e10}
\begin{split}&
\epsilon^2\sum_{j=1}^{l_1} \frac{\sigma_j^2}{1 + 2\alpha P_0^{-1}\epsilon^2\sigma_j^2 j^{2\alpha+1}} \asymp \epsilon^2\sum_{j=1}^{l_1} \sigma_j^2 \\& <C\delta_1\epsilon^2\sum_{j=l_1}^{k_\epsilon}\sigma_j^2 \asymp C\delta_1 \epsilon^2\sum_{j=l_1}^{k_\epsilon} \frac{\sigma_j^2}{1 + 2\alpha P_0^{-1}\epsilon^2\sigma_j^2 j^{2\alpha+1}}
\end{split}
\end{equation}
We have
\begin{equation}\label{e11}
\begin{split}&
\epsilon^2\sum_{j=l_2}^\infty\,\frac{\sigma_j^2}{1 + (2\alpha P_0)^{-1}\,\epsilon^2\,\sigma_j^2 j^{2r+1}} \asymp \epsilon^2\,\sum_{j=l_2}^\infty\,j^{-2\alpha-1}\\& \le
\epsilon^2\,\delta_1^{2\alpha}\,C \,\sum_{k_\epsilon}^{l_2}\,j^{-2\alpha-1} \asymp \epsilon^2\,\delta_1^{2\alpha}\,C\,\sum_{k_\epsilon}^{l_2} \frac{\sigma_j^2}{1 + (2\alpha P_0)^{-1}\,\epsilon^2\,\sigma_j^2\,j^{2\alpha+1}}.
\end{split}
\end{equation}
Now (\ref{e9})-(\ref{e11}) imply that $R_\epsilon - I(P_0-\delta) \to 0$ for some $\delta = \delta(\epsilon) \to 0$  and $\delta_1 = \delta_1(\epsilon) \to 0$ as $\epsilon \to 0$.
\subsection{Proof of Theorem \ref{t3}}
The reasoning is based on the following Lemma.
\begin{lemma}\label{l31}
\begin{equation}\label{st1}
\sup_{x \in B(\alpha,P_0)} E_x ||\tilde x_\epsilon -x||^2 = E_{\theta_\epsilon}||\tilde x_\epsilon - \theta_\epsilon||^2
\end{equation}
with $\theta_\epsilon = \{\theta_{\epsilon k}\}_{k=1}^\infty,  \theta_{\epsilon k} = P_0(a_k - a_{k+1})$.
\end{lemma}
{\sl Proof of Lemma \ref{l31}}. Denote $u_k = a_k^{-1}\sum_{j=k}^\infty \theta_j^2.$ Then
$$
\theta_k^2 = a_k u_k - a_{k+1}u_{k+1}.
$$
Denote $l = [\mu^{-\frac{1}{\beta}}]$.

 We have
\begin{equation}\label{st2}
 E_x ||\tilde x_\epsilon -x||^2 = \mu^2 \sum_{j=1}^l b_j^2 x_j^2 + \sum_{j=l+1}^\infty x_j^2 + \epsilon^2\sum_{j=1}^l \lambda_j^2 \doteq J_1  + J_2 +J_3
 \end{equation}
 respectively.

 We have
 \begin{equation}\label{st3}
 \begin{split}&
 J_1 + J_2 = \mu^2 \sum_{j=1}^l b_j^2 (a_j u_j - a_{j+1} u_{j+1}) + a_{l+1} u_{l+1}\\& =
 \mu^2 a_1 b_1^2 u_1^2 - \mu^2 a_{l+1} b_{l}^2 u_{l+1}^2 \\&
 +   \mu^2 \sum_{j=2}^l a_j u_j(b_j^2 -b_{j-1}^2) + a_{l+1} u_{l+1}.
 \end{split}
 \end{equation}
 The maximum of right-hand side of (\ref{st3}) is attained  for $u_j =  P_0$, $1 \le j < \infty$ with $x_j^2 = P_0(a_j - a_{j+1})$.

 By straightforward calculations, we get
 $J_3= C\epsilon^2 l$.

 If $\beta >\alpha$, we get
 $$
 J_1 + J_2 = \frac{\beta}{\beta-\alpha} l^{-2\alpha}(1 + o(1)).
 $$
 If $\alpha > \beta$, we get
 $$
 J_1 + J_2 = P_0l^{-2\beta} C_1(1 +o(1)).
 $$
 If $\alpha= \beta$, we get
 $$
 J_1 + J_2 = \alpha P_0 l^{-2\alpha} \ln l
 $$
 Minimizing $ J_1 + J_2 + J_3$ with respect  to $l$, we get Theorem \ref{t3}.
 \subsection{Proof of Theorem \ref{t4}}
 It suffices to prove necessary conditions.

 We have
 \begin{equation}\label{}
 \begin{split}&
 E_x ||\tilde x_{\epsilon\mu} - x ||^2 = \epsilon^2 \sum_{j=1}^l (1 - l^{-\beta}j^\beta) \\&
 + l^{-2\beta} \sum_{j=1}^l j^{2\beta}x_j^2 + \sum_{j=l}^\infty x_j^2 \\&
 \ge C\epsilon^2 l + l^{-2\beta} \sum_{j=1}^l j^{2\beta} x_j^2 \doteq J_\epsilon(l,x).
 \end{split}
 \end{equation}
 It easy to seå that, if
 \begin{equation}\label{}
 \sum_{j=1}^l j^{2\beta}x_j^2 \to \infty \quad \mbox{as} \quad l \to \infty,
 \end{equation}
 then
 \begin{equation}\label{}
 \lim_{\epsilon \to 0} \epsilon^{-\frac{4\beta}{1 + 2\beta}} \inf_l J_\epsilon(l,x) = \infty
 \end{equation}

\end{document}